\documentclass[12pt]{extarticle}
\usepackage{mathptmx}

\newcommand{\ind}{\mbox{\mbox{ } \mbox{ } \mbox{ }}}

\newcommand{\be}{\begin{enumerate}}
\newcommand{\ee}{\end{enumerate}}
\newcommand{\bc}{\begin{center}}
\newcommand{\ec}{\end{center}}
\newcommand{\beq}{\begin{equation}}
\newcommand{\eeq}{\end{equation}}
\newcommand{\bqn}{\begin{eqnarray}}
\newcommand{\eqn}{\end{eqnarray}}
\newcommand{\bqns}{\begin{eqnarray*}}
\newcommand{\eqns}{\end{eqnarray*}}
\newcommand{\N}{{\bf N}}
\newcommand{\R}{{\bf R}}

\newcommand{\rarrow}{\rightarrow}

\newcommand{\se}{{\cal E}}

\newcommand{\bdoc}{\begin{document}}
\newcommand{\edoc}{\end{document}}
\newcommand{\docart}{\documentstyle[12pt]{article}}
\newcommand{\vs}{\vspace{.16in}}
\newcommand{\half}{\frac{1}{2}}
\newcommand{\diffeo}{\mbox{ diffeomorphism }} 
\newcommand{\dnorm}[1]{\mbox{$\mid \mid #1 \mid \mid$}}
\newcommand{\twoover}[2]{{\footnotesize \begin{array}[t]{c}\mbox{ {\normalsize #1}}  \\
                               #2 
        \end{array} } \hspace{1cm} }
\newcommand{\threeover}[3]{ {\footnotesize \begin{array}[t]{c}\mbox{ {\normalsize #1}}  \\
                               #2  \\
                                 #3
        \end{array} } \hspace{1cm} }
\newcommand{\fourover}[4]{ {\footnotesize \begin{array}[t]{c}\mbox{ {\normalsize #1}}  \\
                               #2  \\
                                 #3 \\
				#4
        \end{array} } \hspace{1cm}  }

\newenvironment{vover3}[3]{   
       {\footnotesize \begin{array}[t]{c}
                               {\normalsize #1} \vspace{.2cm} \\ 
                                 #2 \\
                                 #3  
                       \end{array} }
                            } { \hspace{1cm} } 
\newenvironment{vover4}[4] {
       {\footnotesize \begin{array}[t]{c} 
                               {\normalsize #1} \vspace{.2cm} \\  
                                 #2 \\
                                 #3 \\
                                 #4  
                       \end{array} }
                            } {  \hspace{1cm} }

\newenvironment{seteqover}[3]{#1 \ \ = \hspace{.6cm}
                                #2 #3} {\ind} 

\newcommand{\norm}[1]{\mbox{$\mid #1 \mid$}}
\newcommand{\Norm}[1]{\mbox{$\parallel #1 \parallel$}}
\newcommand{\imp}{\mbox{$\Longrightarrow$}}
\newcommand{\pr}{\mbox{$\prime$}}
\newcommand{\br}{\mbox{\bf R}}
\newcommand{\fact}{\mbox{$\verb+!+$}}

\newcommand{\tag}[2]{\vspace{.15in}$(#1) \hfill #2 \hfill $ \vspace{.15in} }

\newcommand{\ba}{\mbox{${\bf a}$}}
\newcommand{\bb}{\mbox{${\bf b}$}}
\newcommand{\bi}{\mbox{${\bf i}$}}
\newcommand{\bj}{\mbox{${\bf j}$}}
\newcommand{\bfi}{\mbox{${\bf i}$}}
\newcommand{\bfj}{\mbox{${\bf j}$}}
\newcommand{\bfx}{\mbox{${\bf x}$}}
\newcommand{\bfy}{\mbox{${\bf y}$}}
\newcommand{\bfg}{\mbox{${\bf g}$}}
\newcommand{\bfc}{\mbox{${\bf c}$}}
\newcommand{\bS}{\mbox{${\bf S}$}}
\newcommand{\bV}{\mbox{${\bf V}$}}
\newcommand{\fh}{\hat{f}}
\newcommand{\Ih}{\hat{I}}
\newcommand{\Lh}{\hat{L}}

\newcommand{\alp}{\mbox{$\alpha$}}
\newcommand{\gra}{\mbox{$\alpha$}}

\newcommand{\bet}{\mbox{$\beta$}}
\newcommand{\grb}{\mbox{$\beta$}}

\newcommand{\gam}{\mbox{$\gamma$}}
\newcommand{\Gam}{\mbox{$\Gamma$}}
\newcommand{\del}{\mbox{$\delta$}}
\newcommand{\Del}{\mbox{$\Delta$}}
\newcommand{\sig}{\mbox{$\sigma$}}
\newcommand{\Sig}{\mbox{$\Sigma$}}

\newcommand{\grg}{\mbox{$\gamma$}}
\newcommand{\Grg}{\mbox{$\Gamma$}}
\newcommand{\grd}{\mbox{$\delta$}}
\newcommand{\Grd}{\mbox{$\Delta$}}
\newcommand{\grs}{\mbox{$\sigma$}}
\newcommand{\Grs}{\mbox{$\Sigma$}}
\newcommand{\grf}{\mbox{$\phi$}}
\newcommand{\Grf}{\mbox{$\Phi$}}
\newcommand{\grth}{\mbox{$\theta$}}
\newcommand{\Grth}{\mbox{$\Theta$}}
\newcommand{\gro}{\mbox{$\omega$}}
\newcommand{\Gro}{\mbox{$\Omega$}}
\newcommand{\gre}{\mbox{$\epsilon$}}
\newcommand{\grve}{\mbox{$\varepsilon$}}
\newcommand{\grr}{\mbox{$\rho$}}
\newcommand{\Grr}{\mbox{$\Rho$}}
\newcommand{\grt}{\mbox{$\tau$}}
\newcommand{\gri}{\mbox{$\iota$}}
\newcommand{\grp}{\mbox{$\pi$}}
\newcommand{\grP}{\mbox{$\Pi$}}
\newcommand{\grl}{\mbox{$\lambda$}}
\newcommand{\grL}{\mbox{$\Lambda$}}
\newcommand{\Grl}{\mbox{$\Lambda$}}
\newcommand{\grz}{\mbox{$\zeta$}}

\newcommand{\thet}{\mbox{$\theta$}}
\newcommand{\Thet}{\mbox{$\Theta$}}
\newcommand{\ome}{\mbox{$\omega$}}
\newcommand{\Ome}{\mbox{$\Omega$}}
\newcommand{\eps}{\mbox{$\epsilon$}}

\newcommand{\iot}{\mbox{$\iota$}}

\newcommand{\lam}{\mbox{$\lambda$}}
\newcommand{\Lam}{\mbox{$\Lambda$}}
\newcommand{\zet}{\mbox{$\zeta$}}

\newcommand{\kap}{\mbox{$\kappa$}}

\newcommand{\grx}{\mbox{$\xi$}}
\newcommand{\Grx}{\mbox{$\Xi$}}
\newcommand{\grc}{\mbox{$\chi$}}
\newcommand{\Grc}{\mbox{$\Chi$}}
\newcommand{\grn}{\mbox{$\nu$}}
\newcommand{\grm}{\mbox{$\mu$}}
\newcommand{\grk}{\mbox{$\kappa$}}

\newtheorem{theorem}{Theorem}[section]
\newtheorem{Cor}[theorem]{Corollary}
\newtheorem{Prop}[theorem]{Proposition}
\newtheorem{Lem}[theorem]{Lemma}
\newtheorem{Rem}[theorem]{Remark}
\newtheorem{Proof}[theorem]{Proof}

\newcommand{\rar}{\rightarrow}

\newcommand{\ta}{\tilde{a}}
\newcommand{\tb}{\tilde{b}}
\newcommand{\tc}{\tilde{c}}
\newcommand{\td}{\tilde{d}}
\newcommand{\te}{\tilde{e}}
\newcommand{\tf}{\tilde{f}}
\newcommand{\tg}{\tilde{g}}
\newcommand{\tih}{\tilde{h}}
\newcommand{\ti}{\tilde{i}}
\newcommand{\tj}{\tilde{j}}
\newcommand{\tk}{\tilde{k}}
\newcommand{\tl}{\tilde{l}}
\newcommand{\tm}{\tilde{m}}
\newcommand{\tn}{\tilde{n}}
\newcommand{\tio}{\tilde{o}}
\newcommand{\tp}{\tilde{p}}
\newcommand{\tq}{\tilde{q}}
\newcommand{\tr}{\tilde{r}}
\newcommand{\ts}{\tilde{s}}
\newcommand{\tit}{\tilde{t}}
\newcommand{\tu}{\tilde{u}}
\newcommand{\tv}{\tilde{v}}
\newcommand{\tw}{\tilde{w}}
\newcommand{\tx}{\tilde{x}}
\newcommand{\ty}{\tilde{y}}
\newcommand{\tz}{\tilde{z}}
\newcommand{\tA}{\tilde{A}}
\newcommand{\tB}{\tilde{B}}
\newcommand{\tC}{\tilde{C}}
\newcommand{\tD}{\tilde{D}}
\newcommand{\tE}{\tilde{E}}
\newcommand{\tF}{\tilde{F}}
\newcommand{\tG}{\tilde{G}}
\newcommand{\tH}{\tilde{H}}
\newcommand{\tI}{\tilde{I}}
\newcommand{\tJ}{\tilde{J}}
\newcommand{\tK}{\tilde{K}}
\newcommand{\tL}{\tilde{L}}
\newcommand{\tM}{\tilde{M}}
\newcommand{\tN}{\tilde{N}}
\newcommand{\tO}{\tilde{O}}
\newcommand{\tP}{\tilde{P}}
\newcommand{\tQ}{\tilde{Q}}
\newcommand{\tR}{\tilde{R}}
\newcommand{\tS}{\tilde{S}}
\newcommand{\tT}{\tilde{T}}
\newcommand{\tU}{\tilde{U}}
\newcommand{\tV}{\tilde{V}}
\newcommand{\tW}{\tilde{W}}
\newcommand{\tX}{\tilde{X}}
\newcommand{\tY}{\tilde{Y}}
\newcommand{\tZ}{\tilde{Z}}

\newcommand{\fl}{f_{\grl}}
\newcommand{\el}{E_{\grl}}
\newcommand{\sel}{\se_{\grl}}
\newcommand{\szl}{\sz_{\grl}}
\newcommand{\zl}{Z_{\grl}}
\newcommand{\famfl}{\{ \fl \} }
\newcommand{\sT}{{\cal T}}
\newcommand{\ab}{{\bf a}}
\newcommand{\Def}{{\bf Definition}}
\newcommand{\Pf}{{\bf Proof}}
\newcommand{\uG}{\bar{{\cal G}}}
\newcommand{\sG}{{\cal G}}
\newcommand{\sF}{{\cal F}}
\newcommand{\stG}{\tilde{{\cal G}}}
\newcommand{\tmu}{\tilde{\mu}}
\newcommand{\tpi}{\tilde{\pi}}
\newcommand{\tgrl}{\tilde{\grl}}
\newcommand{\sD}{{\cal D}}
\newcommand{\sDh}{\hat{\sD}}
\newcommand{\sP}{{\cal P}}
\newcommand{\sM}{{\cal M}}
\newcommand{\bsM}{\bar{{\cal M}}}
\newcommand{\sV}{{\cal V}}
\newcommand{\sS}{{\cal S}}
\newcommand{\sI}{{\cal I}}
\newcommand{\sW}{{\cal W}}
\newcommand{\sCZ}{{\cal C}{\cal Z}}
\newcommand{\tsT}{\tilde{\sT}}
\newcommand{\tsD}{\tilde{\sD}}
\newcommand{\tgre}{\tilde{\gre}}
\newcommand{\sH}{{\cal H}}
\newcommand{\tsH}{\tilde{\sH}}
\newcommand{\tsG}{\tilde{\sG}}
\newcommand{\tDel}{\tilde{\Grd}}
\newcommand{\sA}{{\cal A}}
\newcommand{\sB}{{\cal B}}
\newcommand{\sCF}{{\cal CF}}
\newcommand{\bsH}{\bar{\sH}}
\newcommand{\ntwob}{\bar{\nu_2}}
\newcommand{\nb}{\bar{\nu}}
\newcommand{\bF}{\bar{F}}
\newcommand{\bG}{\bar{G}}
\newcommand{\bI}{\bar{I}}
\newcommand{\bZ}{\bar{Z}}
\newcommand{\bh}{\bar{h}}
\newcommand{\bT}{\bar{T}}
\newcommand{\bx}{\bar{x}}
\newcommand{\bxi}{\bar{\xi}}
\newcommand{\bu}{\bar{u}}
\newcommand{\bg}{\bar{g}}
\newcommand{\bN}{\bar{\N}}
\newcommand{\bGrs}{\bar{\Grs}}
\newcommand{\bgrs}{\bar{\grs}}
\newcommand{\bmu}{\bar{\mu}}
\newcommand{\db}{\bar{d}}
\newcommand{\dx}{\dot{x}}
\newcommand{\n}{\norm}
\newcommand{\grfh}{\hat{\grf}}

\newcommand{\Lpfk}{\mbox{$\cap \hspace{-.44cm}\mid$}}
\newcommand{\Emn}{ R_{i_{-m } \ldots i_{-1} , i_0 \ldots i_{n-1}}}
\newcommand{\Rmn}{ R_{i_{-m } \ldots i_{-1} , i_0 \ldots i_{n-1}}}
\newcommand{\En}{ E_{ i_0 \ldots i_{n-1}}}
\newcommand{\Ean}{ E_{ ai_1 \ldots i_{n-1}}}
\newcommand{\Ep}{ E_{ i_0 \ldots i_{p-1}}}
\newcommand{\Sp}{ S_{ i_{0} \ldots i_{p-1}}}
\newcommand{\Sm}{ S_{ {i_{-m } \ldots  i_{-1} }}}
\newcommand{\Sn}{ S_{ i_{0} \ldots i_{n-1}}}
\newcommand{\Gum}{ \Gamma^u_{i_{-m } \ldots  i_{-1} }}
\newcommand{\Gsn}{ \Gamma^s_{ i_0 \ldots i_{n-1}}}
\newcommand{\Em}{ E_{i_{-m } \ldots i_{-1}}}
\newcommand{\Rpn}{ R_{i_{1 } \ldots i_{p} , n}}
\newcommand{\Rkn}{ R_{i_{1 } \ldots i_{k} , n}}
\newcommand{\Ran}{ R_{\bar{\alpha} , n}}
\newcommand{\Enewn}{ E_{ i_1 \ldots i_{n}}}
\newcommand{\Enewp}{ E_{ i_1 \ldots i_{p}}}
\newcommand{\Snewp}{ S_{ i_{1} \ldots i_{p}}}
\newcommand{\Snewn}{ S_{ i_{1} \ldots i_{n}}}
\newcommand{\Snewk}{ S_{ i_{1} \ldots i_{k}}}
\newcommand{\Snewjp}{ S_{ j_{1} \ldots j_{p}}}
\newcommand{\Rnewn}{ R_{i_{1 } \ldots i_{n-1} , i_{n}}}
\newcommand{\Sa}{ S_{\bar{\alpha} }}
\newcommand{\Can}{ C_{\bar{\alpha} , n}}
\newcommand{\Rbm}{ R_{\bar{\beta} , m}}
\newcommand{\Cbm}{ C_{\bar{\beta} , m}}
\newcommand{\Rak}{ R_{\bar{\alpha} , k}}
\newcommand{\fnewn}{ f_{ i_{1} \ldots i_{n}}}
\newcommand{\fnewp}{ f_{ i_{1} \ldots i_{p}}}

\def\bibsameauth{\leavevmode\vrule height .1ex depth 0pt width
2.3em\relax\,}

\begin{document}

\bibliographystyle{plain}

\title{Thermodynamic formalism  for some systems with countable Markov structures}
\author{Michael Jakobson}

% \date{Preliminary Draft---\today} 

\maketitle

  {\centering To the memory of Dmitry Viktorovich Anosov  \par } 
\begin{abstract}
We study ergodic properties of certain piecewise smooth two-dimensional systems
by constructing countable Markov partitions.
Using thermodynamic formalism we prove exponential decay of correleations.
 That extends the results of  \cite{Jak-New-1}, \cite{Jak-New-2}.
Our approach is motivated by
the original method of Anosov and Sinai from \cite{Anosov-Sinai}.

\end{abstract}

\section{Motivation: Folklore Theorem in dimension 1} 

A well-known Folklore Theorem in one-dimensional dynamics
can be formulated as follows. \vs

{\bf Folklore Theorem.} {\it Let $I = [0,1]$ be the
 unit interval, and suppose $\{I_1, I_2, \ldots
\}$ is a  countable collection of disjoint open subintervals of $I$
such that $\bigcup_i I_i$ has the full Lebesgue measure in $I$. 
 Suppose there are constants $K_0 > 1$ and $
 K_1 > 0 $ and mappings $f_i:I_i \rarrow I$ satisfying the
following conditions.

\be
\item $f_i$ extends to a $C^2$ diffeomorphism from the closure of $I_i $ onto
$[0,1]$, and $ \inf_{z \in I_i} \norm{Df_i(z)} > K_0 $ for all $i$. 
\item $ \sup_{z \in I_i}
\frac{\norm{D^2f_i(z)}}{\norm{Df_i(z)}} {\norm{I_i}} < K_1 $ for all $i$.
\ee 
Then, the mapping $F(z)$ defined by $F(z) = f_i(z)$ for $ z \in I_i,$
has a unique invariant  ergodic
probability measure  $\mu$ equivalent to Lebesgue measure on $I$. } \vs

For the proof of the Folklore theorem , the ergodic properties of $\mu$ and the history of the question 
see for example \cite{Adler-79} and \cite{Walters-78}. \\
In  \cite{Jak-New-1} ,  \cite{Jak-New-2} the Folklore Theorem
was generalized to  
two-dimensional maps $F$ which piecewise coincide with certain hyperbolic diffeomorphisms $f_i$.  
 As in the one-dimensional situation there is an essential difference between
a finite and  an infinite number of $f_i$. In the case of an infinite number
 of $f_i$, their derivatives grow with $i$ and  relations between first
and second derivatives become crucial. \\
Models with infinitely many $f_i$ appear when we study
 non-hyperbolic systems, such as quadratic-like maps in dimension 1, and Henon-like maps in 
dimension 2.

\section{Model under consideration. Geometric and hyperbolicity conditions} 
\be
\item
As in \cite{Jak-New-1} ,  \cite{Jak-New-2} we consider the following
$2$-d model. Let $Q$ be the unit square.  Let $\xi = \{E_1, E_2, \ldots, \} $ 
be a countable collection
of  closed curvilinear rectangles in $Q$.  Assume that each $E_i$
lies inside a domain of definition of a $C^2$ diffeomorphism $f_i $
which maps $ E_i $ onto its image $S_i \subset Q$.  We assume each $E_i$
connects the top and the bottom of $Q$. Thus each $E_i$ is bounded from
above and from below by two subintervals of the line segments $\{(x,y) :
y = 1, \ 0 \leq x \leq 1 \}$ and $\{(x,y): y = 0, \ 0 \leq x \leq 1 \}$.
Hyperbolicity conditions that we formulate below imply 
 that the left and right boundaries of $E_i$ are graphs of
smooth functions $x^{(i)}(y)$ with $\left|{dx^{(i)}\over dy}\right| \leq
\gra$ where $\gra$ is a real number satisfying $0 < \gra < 1$. \\
The images $f_i(E_i)=S_i$ are narrow strips connecting the
left and right sides of $Q$ and that they are bounded on the left and
right by the two subintervals of the line segments $\{(x,y): x=0, \ 0
\leq y \leq 1 \}$ and $\{(x,y): x = 1, \ 0 \leq y \leq 1 \}$ and above
and below by the graphs of smooth functions $Y^{i}(X),
\norm{{dY^{(i)}\over dX}} \leq \gra $.  \\
We are saying that $E_i's $ are {\it full height } in $Q$ while the $S_i's $ are
{\it full width} in $Q$. 
\item
For $z \in Q$, let $\ell_z$ be the horizontal line through $z$. We
define $\grd_z(E_i) = diam(\ell_z \bigcap E_i)$, $\grd_{i,max} = \max_{z
\in Q} \grd_z(E_i)$, $\grd_{i,min} = \min_{z \in Q} \grd_z(E_i)$. 
We assume the following \vs

{\bf Geometric conditions}.
\be
\item[G1.]  For $i \neq j$ holds  $int \ E_i \cap int \ E_j =\emptyset$ $int \ S_i \cap int \ S_j =\emptyset$ .
\item[G2.] $mes(Q \setminus \cup_i \  int \, E_i)  = 0$ where $mes$
stands for Lebesgue measure.
\item[G3.] $ - \sum_i \grd_{i,max} \log \grd_{i,min} < \infty$. 
\ee

\item
In the standard coordinate system   for a map 
$F : (x,y) \rightarrow (F_{1}(x,y), F_{2}(x,y))$ we use 
 $ DF(x,y)$ to denote the differential of $F$ at some point
$(x,y) $ and $F_{jx}$, $F_{jy}$, $F_{jxx}$, $F_{jxy}$, 
etc., for partial derivatives of $F_{j},\  j = 1,2$ . 

Let $J_F(z) = \mid F_{1x}(z)F_{2y}(z) - F_{1y}(z)F_{2x}(z) \mid$ be the
absolute value of the Jacobian determinant of $F$ at $z$. \vs

{\bf Hyperbolicity conditions}. 

There exist  constants $0 < \gra < 1$
and $K_0 > 1$ such that for each $i$   the map
$$ 
F(z)  = f_i(z) \mbox{ for } z \in   E_i
$$
satisfies 
\be
\item[H1.] $\n{F_{2x}(z)} + \gra \n{F_{2y}(z)} + \gra^2 \n{F_{1y}(z)}
\leq \gra \n{F_{1x}(z)} $
\item[H2.] $ \n{F_{1x}(z)} - \gra \n{F_{1y}(z)} \geq K_0$.
\item[H3.] $ \n{F_{1y}(z)} + \gra \n{F_{2y}(z)} + \gra^2 \n{F_{2x}(z)}
\leq \gra \n{F_{1x}(z)} $
\item[H4.] $ \n{F_{1x}(z)} - \gra \n{F_{2x}(z)} \geq J_F(z)K_0$.
\ee

 For a real number
$0 < \gra < 1$, we define the cones 

\[  K^u_{\gra} = \{ (v_1, v_2): \norm{v_2} \leq \gra \norm{v_1} \} \]

\[ K^s_{\gra}  = \{ (v_1, v_2): \norm{v_1} \leq \gra \norm{v_2} \} \]

and the corresponding cone fields $K^u_{\gra}(z), K^s_{\gra}(z)$ in
the tangent spaces at points $z \in \R^2$.

The following  proposition proved in \cite {Jak-New-2} relates conditions H1-H4 above with the usual
definition of hyperbolicity in terms of cone conditions. It shows 
that conditions H1 and H2 imply that the $K^u_{\gra}$ cone is mapped into itself by $DF$ and expanded by
a factor no smaller than $K_0$ while H3 and H4 imply that the
$K^s_{\gra}$ cone is mapped into itself by $DF^{-1}$ and expanded by a
factor no smaller than $K_0$. \vs

Unless otherwise stated, we use the $max$ norm on $\R^2$, $\n{(v_1,
v_2)} = \max(\n{v_1}, \n{v_2})$.  

\begin{Prop} \label{cone_cond_prop}
Under conditions H1-H4 above, we have 
\beq \label{ku-cone_1} 
DF(K^u_{\gra}) \subseteq K^u_{\gra}
\eeq

\beq \label{ku-exp_1}
v \in K^u_{\gra} \Rightarrow \n{DFv} \geq  K_0
\n{v}
\eeq

\beq \label{ks-cone_1}
DF^{-1}(K^s_{\gra}) \subseteq K^s_{\gra}
\eeq

\beq \label{ks-exp_1}
v \in K^s_{\gra} \Rightarrow  \n{DF^{-1}v} \geq K_0
\n{v}
\eeq

\end{Prop}

\begin{Rem} \label{Rem0}
{\rm  The first version of hyperbolicity conditions appeared in \cite{Smale-65}. It was developed in particular in
 \cite{Alekseev-1} and
\cite{Hirsch-Pugh} . 
Here we use hyperbolicity conditions from \cite{Jak-New-2}. In \cite{Jak-New-1} we used hyperolicity conditions 
from  \cite{Alekseev-1} which implied the invariance of cones and uniform expansion with
respect to the sum norm $\n{v} = \n{v_1} + \n{v_2}$. }
\end{Rem}
\item
The map 

\[ F(z) = f_i(z) \mbox{ for } z \in int  \ E_i \]

is defined almost everywhere on $Q$.  Let $\tQ_0 = \bigcup_i int \ E_i$,
and, define $\tQ_n, n >0,$ inductively by $\tQ_n = \tQ_0 \bigcap
F^{-1}\tQ_{n-1}$. Let $\tQ = \bigcap_{n \geq 0} \tQ_n$ be the set of
points whose forward orbits always stay in $\bigcup_i int \ E_i$.  Then,
$\tilde{Q}$ has full Lebesgue measure in $Q$, and $F$ maps $\tQ$ into
itself. 

The hyperbolicity conditions H1--H4 imply the estimates on the derivatives of
the boundary curves of $E_i$ and $S_i$ which we described earlier.  They
also imply that any intersection $f_i E_i \bigcap E_j$ is  full width in
$E_j$.  Further, $E_{ij} = E_i \bigcap f_i^{-1}E_j$ is a full height subrectangle
of $E_i$ and $S_{ij} = f_j f_iE_{ij}$ is a full width substrip in $Q$.

Given a finite string $i_0 \ldots i_{n-1}$,  we define inductively

\[ E_{i_0 \ldots i_{n-1}} = E_{i_0} \bigcap f_{i_0}^{-1}E_{i_1i_2 \ldots
i_{n-1}}. \]

Then, each set $E_{i_0 \ldots i_{n-1}}$ 
is a  full height subrectangle of $E_{i_0}$.  

Analogously, for a string
$i_{-m} \ldots i_{-1}$  we define

\[
{\Sm} = f_{ i_{-1}}(S_{i_{-m} \ldots i_{-2}} \bigcap
E_{ i_{-1}}) \] 
 
and get that  $S_{i_{-m} \ldots  i_{-1}}$ is a  full width strip in $Q$.
It is easy to see that $S_{i_{-m} \ldots  i_{-1}} = f_{ i_{-1}} \circ
f_{i_{-2}} \circ \ldots \circ f_{i_{-m}}(E_{i_{-m} \ldots  i_{-1}})$
and that $f_{i_0}^{-1}(S_{i_{-m} \ldots  i_{-1}} )$ is a full-width substrip of
$E_{i_0}$. \\
We also define curvilinear rectangles ${\Rmn}$ by
\[ {\Rmn} = {\Sm} \bigcap {\En} \]
If there are no negative indices then respective rectangle is  full height in $Q$.
For infinite strings, we have the following Proposition. \vs

\begin{Prop} \label{Prop1}
  Any $C^1$  map $F$   satisfying the above geometric
conditions G1--G3 and hyperbolicity conditions H1--H4 has a "topological attractor"
$$ \Lambda =\bigcup_{\ldots i_{-n} \ldots i_{-1}} \bigcap_{k\geq 1} 
S_{i_{-k} \ldots  i_{-1}}$$ \vs

The infinite intersections  $\bigcap_{k=1}^{\infty}S_{i_{-k} \ldots  i_{-1}}$ define  $C^1$    curves
$y(x), \ | dy/dx| \leq  \gra$  which are the unstable manifolds
for the points of the attractor.  The infinite intersections
$\bigcap_{k=1}^{\infty} E_{i_0 \ldots i_{k-1}} $ define  $C^1$
curves 
$x(y), \ | dx/dy| \leq  \gra$  which are the stable manifolds
for the points of the attractor. 
The infinite intersections

 \[ \bigcap_{m=1}^{\infty} \bigcap_{n=1}^{\infty} {\Rmn} \]
define
points of the attractor. 
\end{Prop}
Proposition \ref{Prop1} is a well known fact in hyperbolic theory. For
example it follows from Theorem 1 in \cite{Alekseev-1}. See also
\cite{PS-89}. 
 The union of the stable manifolds 
 has full measure in $Q$. The trajectories of all points
in  this set converge to $\Lambda$. That is the reason to call  $\Lambda$
a topological attractor. 
\item
An $F-$invariant Borel probablility measure
$\mu$ on $Q$ is called a $Sinai-Ruelle-Bowen$ measure (or SRB-measure)
for $F$ if $\mu$ is ergodic and there is a set $A \subset Q$ of
positive Lebesgue measure such that for $x \in A$ and any continuous
real-valued function $\grf:Q \rarrow \R$, we have
\beq \label{time_avg}
 \lim_{n \rarrow \infty} \frac{1}{n} \sum_{k=0}^{n-1} \grf(F^kx) = \int \grf d\mu. 
\eeq
Existence of an SRB measure is a much stronger result, than \ref{Prop1}.
It allows to describe statistical 
properties of trajectories in a set of positive phase volume. It requires some
additional assumptions.
\ee

\section{Distortion conditions } 
 
As we have a countable number of domains the derivatives of $f_i$
grow.
We  formulate certain assumptions on the second
derivatives.   We  use the 
distance function $d((x,y), (x_1,y_1)) = \max(\n{x - x_1}, \n{y -
y_1})$ associated with
the norm $\n{v} =
\max(\n{v_1}, \n{v_2})$ on vectors $v = (v_1, v_2)$.

As above, for a point $z \in Q$, let $l_z$ denote the horizontal line through
$z$, and if $E \subseteq Q $, let $\grd_z(E)$ denote the
diameter of the horizontal section $l_z \bigcap E$. We call $\grd_z(E)$
the $z-width$ of $E$.  

In given coordinate systems we write $f_i(x,y) = (f_{i1}(x,y), f_{i2}(x,y))$. 
We use
$f_{ijx}, f_{ijy}, f_{ijxx}, f_{ijxy}$, etc. for partial derivatives 
of $f_{ij}, j = 1,2$. 

We define 

\[ \n{D^2f_i(z)} = \max_{j=1,2, (k,l) = (x,x), (x,y), (y,y)}
\n{f_{_{ijkl}}(z)}. \]

Next we formulate distortion conditions which
are used to control the fluctuation of the
derivatives of iterates of $F$ along unstable manifolds, and to construct Sinai local measures.
 
Suppose there is a constant $C_0 > 0$ such that the following {\it
distortion conditions} hold

\be
\item[D1.] $ \ \  \sup_{z \in E_i, i \geq 1}  \frac{\n{D^2f_i(z)}}
 { \norm{f_{i1x}(z)}}\grd_z(E_i) < C_0 $. 
\ee \vs

Our conditions imply the following theorem proved in \cite{Jak-New-1},  \cite{Jak-New-2}.
\begin{theorem} \label{Theorem1}
 Let $F$ be a piecewise smooth mapping as above
satisfying the geometric conditions G1--G3, the hyperbolicity conditions
H1--H4 and  the distortion condition D1. 

Then, $F$ has an $SRB$ measure $\mu$ supported on 
$\Lambda $ whose basin has full Lebesgue measure in
$Q$.   Dynamical system   $(F,\mu)$ satisfies the following properties.\\
\be
\item $(F,\mu)$ is measure-theoretically isomorphic to a 
Bernoulli shift.
\item $F$ has finite entropy with respect to the measure
$\mu$,  and the entropy formula holds
\beq \label{ent_form1}
 h_{\mu}(F) = \int \log|D^uF| d\mu 
\eeq
where $D^uF(z)$ is the norm of the derivative of $F$ in the
unstable direction at $z$. 
\item
\beq \label{ent_form}
 h_{\mu}(F) = \lim_{n \rarrow \infty} \frac{1}{n} \log \n{DF^n(z)} 
\eeq
where the latter limit exists for Lebesgue almost all $z$ and is
independent of such $z$. 

\ee
\end{theorem}

\section{Additional hyperbolicity and distortion conditions and statement of the main  theorem}
When applying  thermodynamic formalism to hyperbolic attractors one considers the function 
 $\phi(z) = - \log(D^u F(z)) $.
Thermodynamic formalism is based on the fact that the pullback of  $\phi(z)$
 into a   symbolic space determined by some Markov
partition is a  locally H{\" o}lder function. \\ 
We prove H\"older property of $\phi(z)$ assuming an extra  hyperbolicity condition, 
and a distortion condition D2 stronger than D1. 

{\bf Hyperbolicity  condition  H5}. \vs
\be
\item[H5.] $ \ \   \frac{1}{K_0^2} + \gra^2 < 1$. 
\ee \vs

{\bf Distortion condition  D2}. \vs
\be
\item[D2.] $ \ \  \sup_{z \in E_i, i \geq 1}  \frac{\n{D^2f_i(z)}}
 { \norm{f_{i1x}(z)}} < C_0 $. 
\ee \vs

\begin{Rem} Condition D2 is too strong to be useful for quadratic-like systems.
In dimension 1 it reads as $\n{\frac{F_{ixx}}{F_{ix}}}< c$ instead of 
$\n{\frac{F_{ixx}}{F^2_{ix}}}< c$. However instead of D2 one can assume 
additional hyperbolicity conditions, which can be vaguely formulated as
"contraction of $f_i$ grows faster than expansion" . That approach will be discussed in a forthcoming paper.
\end{Rem}
Assuming additionally H5 and D2 we prove that H\"older  functions have exponential decay of correlations. \\
Let ${\cal H}_{\gamma}$ be the space of functions on $Q$ satisfying 
H\"older property with exponent $\gamma$ 
\[ \norm{\phi(x)-\phi(y)} \le c\norm{x-y}^{\gamma}\]
Then the following theorem holds. \\

\begin{theorem} \label{Main}
 Let $F$ be a piecewise smooth mapping as above
satisfying the geometric conditions G1--G3, the hyperbolicity conditions
H1--H5 and  the distortion condition D2. 
 Then
 $(F,\mu)$ has exponential decay of correlations  for $\phi, \psi \in {\cal H}_{\gamma}$.
Namely there exist $\eta(\gamma)<1$ and
  $C = C(\phi, \psi )$
such that
\beq \label{decay}
 \norm{\int \phi(\psi \circ F^n) d\mu - \int \phi d\mu \int \psi d\mu} < C\eta^n
\eeq
\end{theorem}

\section{H{\" o}lder  properties of  $\log(D^u F(z)) $ }
\be
\item

Although Markov partitions are partitions of the attractor, we need to check  H{\" o}lder property
on actual two-dimensional  curvilinear rectangles 
${\Emn}$. We  call respective partition Markov as well.
In our  model Markov partition consists of initial full height rectangles $E_i$. \\
We consider rectangles  ${\Emn}$ with  $m \ge 0, n \ge 1$. 
We use notation $m=0$ if there are no negative coordinates, which means ${\Emn}= R_{i_0, \ldots ,i_{n-1}}$
is a full height rectangle. The assumption  $n \ge 1$,  means that variations are measured 
between points which belong to the same full height rectangle. \\
By definition the function $\log  D^uF$ is {\bf locally  H{\" o}lder} if for $m \ge 0$, $n \ge 1$
  the variation of $\log  D^uF$ on  ${\Emn}$ satisfies
\beq \label {var1}
 var(\log D^uF) | {\Emn}< C\theta_0^{min(m,n)}
\eeq
for some $C>0$, $\theta_0 < 1$.\\
\begin{Prop} \label{Hproposition}
 $\log  D^uF$ is a locally  H{\" o}lder function.
\end{Prop}
We prove Proposition \ref{Hproposition} with some $\theta_0 $ and  $C$ determined by hyperbolicity and distortion conditions.
\be
\item
The sets ${\Emn} $ are bounded from above and below by 
some arcs of two unstable curves ${\Gum}$, which are
images of some pieces of the top and bottom of $\tilde Q$,
and  from left and right by some arcs
 of two  stable curves ${\Gsn}$, which are  preimages of some pieces the left and right boundaries of $\tilde Q$ . \\
Let $Z_1, Z_2 \in {\Rmn}$ be two 
points on the attractor. We connect $Z_1,Z_2$ by two pieces of their unstable manifolds to two points $Z_3,Z_4$ which belong to the same stable manifold. Let

$\gamma_1 = \gamma(Z_1,Z_3) \subset W^u(Z_1)$, $\gamma_2 = \gamma(Z_2,Z_4) \subset W^u(Z_2)$,
$\gamma_3 = \gamma(Z_3,Z_4) \subset W^s(Z_3)$ be respective curves all located inside
 ${\Rmn}$.\\
 We estimate
\begin{eqnarray*} \label {connectz1z2}
\norm{ \log D^uF(Z_1)-  \log D^uF(Z_2)} \le \norm{\log  D^uF(Z_1)- \log D^uF(Z_3)} +\\
 \norm{ \log D^uF(Z_3)- \log  D^uF(Z_4)}+
\norm{\log  D^uF(Z_4)- \log  D^uF(Z_2)} 
\end{eqnarray*}

\item First we estimate $\norm{\log  D^uF(Z_1)- \log D^uF(Z_3)}  $. 
 We connect $Z_1$ and $Z_3$
by a chain of small rectangles $R \subset {\Rmn}$ covering  $\gamma_1$.
Then $\norm{\log  D^uF(Z_1)- \log D^uF(Z_3)}  $ is majorated by the sum of
similar differences for points  $z_1, z_2 \in W^u(Z_1) \subset R$. \\
Because of cone conditions we can choose rectangles $R= \Delta x \times \Delta y$
satisfying $\norm{\Delta y} < \gra \norm{\Delta x}$.
 Let $R$ be one of such rectangles.\\
%%%%%%%%%%%%%%%%%%%
%%%%%%%%%%%%%%%%%%%%%%%%
 Hyperbolicity conditions imply the  following properties, see \cite{Jak-New-2}.
\be
\item Any unit  vector in $K^u_{\gra}$
at a point $z \in E_i$,  in particular a tangent vector to $W^u(z)$,
has coordinates $(1,a_z)$ with $\n{a_z} < \gra$. 
\item
\beq \label {HJN1}
\n{D^uF(z)} = \n {F_{1x}(z)+a_zF_{1y}(z)}
\eeq
\item
\beq \label {Hyp1}
\frac{\n{F_{1y}}}{\n{F_{1x}}} < \gra 
\eeq
\item
\beq \label {Hyp2}
\frac{\n{F_{2x}}} {\n{F_{1x}}} < \gra 
\eeq
\item
\beq \label {Hyp3}
\frac{\n{F_{2y}}}{\n{F_{1x}}} < \frac{1}{K_0^2} + \gra^2 
\eeq
\ee

Assuming without loss of generality $F_{1x} >0$  for all $x \in E_i$ we get
that variation of $\log\n{D^uF}$ between two  points $z_1, z_2 \in W^u(Z_1) \subset R$
equals
\beq \label {log0}
 \log \big [F_{1x}(z_1) \big (1+a_{z_1}\frac{F_{1y}}{F_{1x}}(z_1)\big ) \big ]-
\log \big [F_{1x}(z_2) \big (1+a_{z_2}\frac{F_{1y}}{F_{1x}}(z_2)\big )\big ] 
\eeq
We split it into two expressions and estimate separately
\beq \label {log1}
 \log F_{1x}(z_1)-  \log F_{1x}(z_2)
\eeq
and
\beq \label {log2}
\log \big (1+a_{z_1}\frac{F_{1y}}{F_{1x}}(z_1)\big ) -
\log \big (1+a_{z_2}\frac{F_{1y}}{F_{1x}}(z_2)\big )
\eeq
We rewrite \ref{log2} as
\beq \label {log3}
\log  \big (1+ \frac { a_{z_1}\frac{F_{1y}}{F_{1x}}(z_1) -   a_{z_2}\frac{F_{1y}}{F_{1x}}(z_2) } 
{1+a_{z_2}\frac{F_{1y}}{F_{1x}}(z_2)} \big )
\eeq
As denominator of the fraction in \ref{log3} is uniformly bounded away from 0, we stimate the numerator and rewrite it as a sum of two expressions
\beq \label {logar}
\n{a_{z_1}}\frac {\n{F_{1y}(z_1)F_{1x}(z_2)- F_{1y}(z_2)F_{1x}(z_1)}}{F_{1x}(z_1)F_{1x}(z_2)}
\eeq
and 
\beq \label {log4}
\n{a_{z_1}- a_{z_2}}\n{ \frac {F_{1y}(z_2)}{F_{1x}(z_2)}}
\eeq
As $C^2$ sizes of unstable manifolds are uniformly bounded ( see  \cite{Jak-New-2}),
\ref{log4} is bounded by $c\n{\Delta x}$.
We rewrite \ref{logar} as 

\beq \label {log5}
a_{z_1} \big [  
\frac {F_{1y}(z_1)\big (F_{1x}(z_2)-F_{1x}(z_1) \big )}{F_{1x}(z_1)F_{1x}(z_2)} 
+ \frac { F_{1y}(z_1)- F_{1y}(z_2)}{F_{1x}(z_2)} \big ]
\eeq
As $\n{\frac{F_{1y}(z_1)}{F_{1x}(z_1)} } < \gra$, both expressions are estimated similarly.\\
 As we are moving along 
$W^u$, we get  $\n{\Delta y}<\gra \n{\Delta x}$. \\
We use the mean value theorem and  
 distortion assumptions, and get estimates bounded by
\beq \label {log6} 
c\n{\Delta x}\frac{F_{1x}(\theta)}{F_{1x}(z_2)}
\eeq
Then it remains to estimate $\frac{F_{1x}(\theta)}{F_{1x}(z_2)}$ or equivalently
$ \n{\log F_{1x}(\theta) - \log F_{1x}(z_2)} $,
which is the same estimate as \ref{log1}.\\
In order to estimate \ref{log1}
 we use again the mean value theorem 
and  distortion assumptions. \\
Then we get
\beq \label {log7}
\n{\log F_{1x}(\theta) - \log F_{1x}(z_2)} < c \n{\Delta x}
\eeq
and respectively
\beq \label {log8}
\frac{F_{1x}(\theta)}{F_{1x}(z_2)}< \exp ( c \n{\Delta x})
\eeq
We combine the previous estimates and get 
\beq \label {log8} 
\norm{\log  D^uF(z_1)- \log D^uF(z_3)}<  c\n{\gamma(z_1,z_3)}
\eeq
From hyperbolicity conditions we get
\beq \label{unstablelegth}
\n{\gamma(z_1,z_3)} < C_2 \frac{1}{K_0^n}
\eeq
That implies
\beq \label{stableconnection1}
\norm{\log  D^uF(z_1)- \log D^uF(z_3)} < C_2 \frac{1}{K_0^n}
\eeq
where $C_2$ is a uniform constant.\\
Similar inequality holds for $\gamma(z_2,z_4)$. \\
\beq \label{stableconnection2}
\norm{ \log D^uF(z_2)- \log D^uF(z_4)} < C_2 \frac{1}{K_0^n}
\eeq

\item
Next we  estimate the variation of $\log  \n{D^uF(z)} $
 between points $Z_3$ and $Z_4$, which belong to the same stable manifold
$W^s(Z_3)=W^s(Z_4) \subset {\Emn}$. Thus we need to estimate
\beq \label {HJN2}
\log \n {F_{1x}(Z_3)+a_{Z_3}F_{1y}(Z_3)} - \log \n{ F_{1x}(Z_4)+a_{Z_4}F_{1y}(Z_4)}
\eeq
As above we split the variation \ref{HJN2} into \ref{log1}, \ref{logar} and \ref{log4}.\\
This time 
instead of moving
along $W^u(Z_1)$ we are moving along  $W^s(Z_3)$, which connects $Z_3$ and $Z_4$.
In that case we use $\n{\Delta x}< \gra \n{\Delta y}$, 
so $\Delta y$ variations are added. As above estimates \ref{log1}, \ref{logar} contribute
less than  
\beq \label{stableconnection4}
c \n{\gamma_3} < C_2 \frac{1}{K_0^m}
\eeq
Note that the Lipschitz dependence 
of the unstable leaves from $y \in W^s_0$ is
unclear in our setting,
differently from the classical case in dimension two. \\
The following  lemma related to that issue is sufficient for our purposes.
\begin{Lem} \label{angles}
There exist $c_0 >0$, $0<\theta_0 < 1$ such that
\beq \label{angles1}
\n{a_{Z_3}-a_{Z_4}} < c_0 \theta_0^m
\eeq
\end{Lem}
Proof. \\
 
We assume by induction that for any rectangle  ${\Rmn}$, and for any points $Z_3,Z_4 \in {\Rmn}$ 
 of intersection of two unstable manifolds
$W^u_1$, $W^u_2$ with the same stable manifold $W^s_0$,
the inequality  \ref{angles1}
holds.  Then we prove 
\beq \label{angles2}
\n{a_{F(Z_3)}-a_{F(Z_4)}} < c_0 \theta_0^{m+1}
\eeq
$DF$ maps a unit vector $\vec{v}=(1,a)$ into 
$(F_{1x}+F_{1y}a, F_{2x}+F_{2y}a)$. Then the normalized vector $DF\vec{v}$ has
second coordinate
\beq \label {angles3}
a' = \frac{ \frac{F_{2x}}{F_{1x}}+ \frac{F_{2y}}{F_{1x}}a}{1+ \frac{F_{1y}}{F_{1x}}a}
\eeq
We denote $Z_3 = z$, $Z_4=w$ and estimate 
\beq \label {angles4}
 \frac{ \frac{F_{2x}}{F_{1x}}(z)+ \frac{F_{2y}}{F_{1x}}(z)a(z)}{1+ \frac{F_{1y}}{F_{1x}}(z)a(z)}
-  \frac{ \frac{F_{2x}}{F_{1x}}(w)+ \frac{F_{2y}}{F_{1x}}(w)a(w)}{1+ \frac{F_{1y}}{F_{1x}}(w)a(w)}
\eeq
After cross multiplying we get denominator  bounded away from $0$. Therefore it is enough
  to estimate two terms 
\beq \label {angles5}
  \frac{F_{2x}}{F_{1x}}(w)\big ({1+ \frac{F_{1y}}{F_{1x}}(z)a(z)}\big ) -
 \frac{F_{2x}}{F_{1x}}(z)\big ({1+ \frac{F_{1y}}{F_{1x}}(w)a(w)}\big ) 
\eeq
and
\beq \label {angles6}
 \frac{F_{2y}}{F_{1x}}(w)a(w)\big (1+ \frac{F_{1y}}{F_{1x}}(z)a(z)\big ) -
 \frac{F_{2y}}{F_{1x}}(z)a(z)\big (1+ \frac{F_{1y}}{F_{1x}}(w)a(w)\big )
\eeq
Both expressions are estimated similarly. To estimate \ref{angles6} we split it into
\beq \label {angles7}
 \frac{F_{2y}}{F_{1x}}(w)a(w)-\frac{F_{2y}}{F_{1x}}(z)a(z)
\eeq
and
\beq \label {angles8}
a(z)a(w)\big(  \frac{F_{2y}}{F_{1x}}(w) \frac{F_{1y}}{F_{1x}}(z)-
\frac{F_{2y}}{F_{1x}}(z) \frac{F_{1y}}{F_{1x}}(w)\big )
\eeq
As above we use elementary algebra and get expressions of the type
\beq \label {angles7a}
\frac{F_{1x}(w) - F_{1x}(z)}{F_{1x}(z)}
\eeq

and 
\beq \label {angles7b}
\frac{F_{2y}(w) - F_{2y}(z)}{F_{1x}(z)}
\eeq

We split $\gamma_3$ into small intervals, and apply the mean value theorem.
The ratios $ \frac{F_{1x}(\theta) }{F_{1x}(z)} $ or equivalently the differences
$ \log F_{1x}(\theta) - \log F_{1x}(z) $ for close points $\theta, z$ on the same
stable manifold are estimated (using again the mean value theorem and D2)
as
\beq \label {angles7c}
\log F_{1x}(\theta) - \log F_{1x}(z) < C_0(1+\gra)\Delta y
\eeq
Thus for any two points $z$ and $\theta$ on the same stable manifold
\beq \label {angles7d}
\log F_{1x}(\theta) - \log F_{1x}(z) < C\n{z-\theta}
\eeq
In particular for all points $z$ and $\theta$ on the same stable manifold
the ratios $\frac{F_{1x}(\theta)}{F_{1x}(z)}$ are uniformly bounded.\\
Thus estimate \ref {angles8} contributes
\beq \label {angles9}
 C \n{\gamma_3}
\eeq
When estimating \ref{angles7} we get similar terms estimated as \ref{angles9}, and 
\beq \label {angles10}
 \frac{F_{2y}}{F_{1x}}(z)\big ( a(z)-a(w) \big )
\eeq
After we combine all terms except  \ref{angles10} we get an estimate 
\beq \label {angles11}
M_0C_0 \frac{1}{K_0^m}
\eeq
where $M_0$ is a uniform constant, which depends on the number of similar terms
that we added above, and $C_0$ is the distortion constant from condition
D2. For \ref{angles10}  we use inductive assumption \ref{angles1} and
get  a total estimate
\beq \label{angles12}
\n{a_{F(Z_3)}-a_{F(Z_4)}} < M_0C_0 \frac{1}{K_0^m} + \big (\frac{1}{K_0^2}+\gra^2\big)c_0 \theta_0^m
\eeq
As $K_0 > 1$ we can choose $\theta_0 < 1$ satisfying for some $A_0 > 1$
\beq \label {angles13}
 \theta_0 = \frac{A_0}{K_0}
\eeq
Also H5 implies that we can choose  $\theta_0 < 1$ satisfying simulteneously

\beq \label {angles13}
\frac{1}{K_0^2}+\gra^2  < \theta_0 
\eeq
Then if
\beq \label {angles14}
c_0 >\frac{ M_0C_0}{\theta_0 - \frac{1}{K_0^2}+\gra^2  }
\eeq
we get the left side of \ref{angles12}  less than $c_0 \theta_0^{m+1} $. \\
Q.E.D. \\
From Lemma \ref{angles} and \ref{stableconnection4} we get

\beq \label{stableconnection3}
\norm { \log D^uF(z_3)- \log D^uF(z_4)}  < C_3 \theta_0^m
\eeq

Combining  \ref{stableconnection1},  \ref{stableconnection2},   \ref{stableconnection3}
  we conclude the proof of Proposition \ref{Hproposition}.
 \ee
\item
We combine several corollaries from Proposition \ref{Hproposition} and from the arguments
used in its proof. 
\begin{Cor} \label{c1}
There exists  $c$ independent of $i$ such that for any $z_1,z_2 \in E_i $ holds  
\beq  \label {bndratiof1x}
\frac{\norm{ f_{i1x}(z_1)} } { \norm{ f_{i1x}(z_2)}} < c
\eeq
\end{Cor}
Here  $z_1,z_2$ do not need to be on the attractor. \\
 To prove Corollary \ref{c1} we fix an arbitrary stable manifold
$W^s_i \subset E_i$, and connect $z_1$ to $z_3 \in W^s_i$ by a horizontal segment $\sigma$.
As full height rectangles are bounded from above and below by horizontal segments,
$\sigma$ lies entirely  in $E_i$. Similarly we connect 
 $z_2$ to $z_4 \in W^s_i$. Then  \ref{log7} and \ref{angles7d} imply \ref{bndratiof1x}.\\
Let $\delta_z(E_i)$ be the width of the horizontal crossection of $E_i$ through $z \in E_i \cap \Lambda$. As $\delta_z(E_i)$ are mapped onto full width unstable curves 
we get  from \ref {bndratiof1x} 
\begin{Cor}
There exists  $c$ independent of $i$ such that for any $z_1,z_2 \in E_i $ holds 
\beq  \label {bndratiowidths2}
\frac{\delta_{z_1}(E_i) }{\delta_{z_2}(E_i) } < c
\eeq
\end{Cor}
\begin{Rem}
Property \ref{bndratiowidths2} demonstrates restrictions on geometry imposed by condition D2.\\
Conditions of Theorem \ref{Theorem1} allow widths of $E_i$ to oscillate
exponentially between $a^i$ and $b^i$ for some $0<a<b<1$. However  from \ref{bndratiowidths2}
 we get that
ratios are uniformly bounded.
\end{Rem}
Applying \ref{var1} to the full height  rectangles $E_i$ we get
for all $z_1,z_2 \in E_i  \cap \Lambda$   
\begin{Cor}

\beq  \label{varEi} 
 var(\log D^uF) | E_i < C
\eeq
\end{Cor}
and
\begin{Cor}
\beq  \label {bndratioDuf}
\frac{D^u f_{i}(z_1)}  { D^u f_{i}(z_2)} < c
\eeq
\end{Cor}
For any $z$ on an unstable curve  $W^u(z) \subset  E_i$
which is full width in $E_i$   let $\norm{W^u(z,E_i)}$  be the respective length .
As   $\norm{W^u(z,E_i)}$ coincide up to a uniformly bounded factor
with $ \frac{1}{D^u f_i(z)}$
 we get from  \ref {bndratioDuf} 
\begin{Cor}
There exists  $c$ independent of $i$ such that for any $z_1,z_2 \in E_i\cap \Lambda$ holds  
\beq  \label {bndratiowidths1}
\frac{\norm{ W^u(z_1,E_i)} } { \norm{ W^u(z_2,E_i)}} < c
\eeq
\end{Cor}
Note that \ref {bndratiowidths1} also follows from \ref{bndratiowidths2} because
at a given point $z \in E_i$ ratios : $\frac{\delta_z(E_i)}{\norm{ W^u(z_1,E_i)}}$
are uniformly bounded.

Although the  next statement is not used in the proof of the main theorem ,
it is usefull for understanding  the geometry of partitions into ${\En}$. 
We claim that    \ref{bndratiowidths1} is valid for rectangles of any order. 
\begin{Rem}
{\it Let ${\En}$ be a full height rectangle of order $n$. \\
Then for any two points  $z_1, \ \ z_2 \in {\En} \cap \Lambda$ holds }
{\vspace{.1in}}
\beq  \label {bndratiowidths1n}
\frac{\norm{ W^u(z_1, {\En})} } { \norm{ W^u(z_2, {\En})}} < c
\eeq
\end{Rem}
To prove  \ref {bndratiowidths1n} we rewrite the ratio
$\frac{D^uF^n(z_2)  }{D^uF^n(z_1)  } $ as
\beq \label{productnew}
 \frac {\prod_{i=0}^{n-1}
D^uF(F^iz_2) }{\prod_{i=0}^{n-1} D^uF(F^{i}z_1)} 
\eeq
and consider
\beq \label{sumnew}
\sum_{i=0}^{n-1}  
\norm{\log D^uF(F^iz_2) - \log D^uF(F^{i}z_1) }
\eeq
As in the proof of Proposition \ref{Hproposition} we split the estimate of each term 
into estimates along stable and unstable manifolds in the images 
$F^p({\En}) = R_{i_0 \ldots i_{p-1},i_p \ldots i_{n-1}}$. For each term estimates from the proof of Proposition \ref{Hproposition} imply respective bounds : $C\theta_0^p$ on stable manifolds,
and $C \frac{1}{K_0^{n-p}}$ on unstable manifolds. Thus we get a uniform bound in \ref{sumnew},
which implies \ref{bndratiowidths1n} .

\item
According to \ref{Prop1} points of $\Lambda$ are identified with two-sided sequences
\[ (\ldots i_{-m} \ldots i_{-1},i_0i_1 \ldots i_n \ldots)\] In order to use Ruelle-Bowen approach
  we define a function $\phi^u$ corresponding to 
$- \log D^uF(z)  $ on the space of one-sided sequences. 
We fix some unstable manifold $W^u_0$.
Let $z=(x,y) \in \Lambda$, $z_0=W^s(z)\cap W^u_0$.
 For $ \phi(z) = - \log D^u(z)$
define
\beq \label{series} 
u(z) = \sum_{k=0}^{\infty}\phi(F^k(z))- \phi(F^k(z_0))
\eeq
and
\beq \label{cohomformula}
\psi(z) = \phi(z) - u(z) + u(Fz)
\eeq
From \ref{var1} we get that the series \ref{series} converge uniformly. \\
For $\psi(z)$ all terms with $z$ cancel, so $\psi(z)$ as a function of $z_0$ 
depends only on nonnegative iterates of $z$. \\

Let $\Omega_+ = \{ x= (i_0i_1 \ldots i_n \ldots ) \}$ be the space of one-sided
sequences corresponding to stable manifolds. \\
The above function $\psi(z)$ defined on $\Omega_+$ will be denoted $\phi^u(x)$.
On the space $\Omega_+ $ we use the metric
$d(x^1,x^2) = 2^{-n}$, where $n = min \{ k : i_k^1 \neq i_k^2 \}$. \\
  For a function $\phi(x)$ on the symbolic space $\Omega_+ $ let
\[ V_n(\phi)= \sup{|\phi(x) - \phi(y)| : x_i=y_i, i = 0, \ldots , n-1} \]
$\phi(x)$ is called
locally H{\" o}lder  if there are $C >0, 0 < \theta < 1$ such that $\forall n \ge 1$
\beq \label {varonesided}
 V_n(\phi) < C \theta^n
\eeq
We use the same arguments as in \cite{Bowen} which prove that the H\"older property on the space of two-sided sequences implies the H\"older property for respective function on the space of one-sided sequences. Then Proposition \ref{Hproposition} implies
\begin{Cor} \label{onesidedHolder}
$\phi^u(x)$ is  a locally H{\" o}lder function on the symbolic space $\Omega_+$.
\end{Cor}

\ee

\section {Some sufficient conditions for exponential decay of correlations in countable shifts}
\be
\item
We refer to  \cite{Sarig-1} for the following general results about shifts with countable alphabets
.\\
Let $T$ be the shift transformation on the space $X$ of admissible one-sided sequences
determined by an infinite matrix $A$ with $a_{ij} = 0,1$. Here $i,j$ are states of an infinite
alphabet. Let $C = [i_0, \ldots , i_{n-1}]$ be cylinder sets.
The system $(X,A,T)$ is called {\it topologically mixing} if the following holds.
\beq  \label {topmix}
\forall C_1,C_2 \ \ \exists N(C_1,C_2) : \forall n>N(C_1,C_2) \ \ C_1 \cap T^{-n}C_2 \neq \emptyset
\eeq
 for any two cylinder sets $C_1$ and $C_2$ . \\
Let $(X_A, T)$ be a topologically mixing countable shift, and let $\phi(x)$ be a 
locally H\"older function. \\
Set $\phi_n(x)= \sum_{k=0}^{n-1}\phi \circ T^k(x)$.
Define  $ Z_n(\phi,a) $ by
\beq  \label {statsum}
Z_n(\phi,a) = \sum_{T^n x=x, x_0 = a} e^{\phi_n(x)} 
\eeq
Then
\beq  \label {Potphi}
P(\phi) = \lim_{n \to \infty}\frac{1}{n}\log Z_n(\phi,a) 
\eeq
does not depend of the choice of $a$. \\

{\bf Definition} 
Assuming  $(X_A, \sigma)$ is topologically mixing and $\phi$ is locally  H\"older continuous
$\phi$ is called 
{\bf positive recurrent} if there is  $\lambda > 0$ such that for
 any given symbol $a$ there is a constant $M_a > 1$ and an integer $N_a$
such that for every $n \ge N_a$ holds 
\[ \frac{ Z_n(\phi,a) }{\lambda^n} \in [M_a^{-1}, M_a] \]
Let $L_{\phi}$ be the Ruelle operator
\beq  \label {R}
L_{\phi}f (x) = \sum_{Ty=x} e^{\phi(y)}f(y)
\eeq
The following is a part of Theorem 4 in \cite{Sarig-1}. \\
\begin{theorem} \label{Sarig1} 
Let $(X_A, T)$ be a topologically mixing countable shift, and let $\phi(x)$ be a 
locally H\"older function such that $P(\phi) < \infty$. If
$\phi$ is positive recurrent, then $\lambda= e^{P(\phi)}$   and
there exist a $\sigma$-finite measure $\nu$ and a function $h>0$ such that 
$L^*_{\phi} \nu = \lambda \nu$, \ \ $L_{\phi}h = \lambda h$, \ \ $\nu(h)=1$,
and for every uniformly continuous function $f$ such that $ ||fh^{-1}||_{\infty} < \infty$
holds
\[ \lambda^{-n}L^n_{\phi}f \rar \nu(f)h \] uniformly on compacts.
\end{theorem}
Consider the space of functions $\cal L$ with bounded norm
$ \norm{f}_{\cal L}$ which is the sum of $||f||_{\infty}$ and some fixed H\"older norm.
The next theorem which follows from several results in \cite{Sarig-1} implies
exponential decay of correlations for H\"older functions restricted to the attractor.
\begin{theorem} \label{Sarig2} 
Suppose the following properties are satisfied. \\
a)  $(X_A, T)$ is topologically mixing.\\
b) $P(\phi) < \infty$. \\
c) $\phi$ is {\it positive recurrent}.\\
d) $h$ is bounded away from $0$ and $\infty$. \\

Then there exist $K > 0, \ \ \theta \in (0,1)$ such that for $f \in \cal L$ holds
\beq \label{maindecay}
\norm { \lambda^{-n} L^n_{\phi} f - h \nu (f) }_{\cal L} < K \theta^n \norm{f}_{\cal L} 
\eeq
\end{theorem}

In general for infinite shifts none of the properties a),b),c),d) is automatic. \\
The following  property  was introduced by
Aaronson,Denker,Urbanski (\cite{ADU-1}) and Yuri (\cite{Yuri}).\\
 {\bf Finitely Many Images Property }. \\
The number of different rows of the matrix $A$ is finite. \\

 In the case of  finitely many images  
the following results from \cite{Sarig-1} simplify the set of properties sufficient for
exponential decay of correlations.\\

(i) \ \ If $A$ satisfies  finitely many images property and 
$\norm{L_{\phi}(1)}_{\infty} < \infty$,  then $\phi$ is positive recurrent .\\
(ii) \ \ If $\phi$ is positive recurrent and   $P(\phi) < \infty$, then
$h$ is bounded away from $0$ and infinity.

\item As a corollary from the above results from \cite{Sarig-1} : theorems \ref{Sarig1}, \ref{Sarig2},
and properties (i) and (ii)
 we get
\begin{Prop} \label{sufficient}
Suppose there is a Markov partition of the attractor satisfying the following properties.
\be
\item  The matrix $A$ of admissible transitions is topologically mixing and has finitely many different rows. 
\item  $\Phi(x,y) = - log \norm {D^uF}$
is H\"older on  the space of admissible sequences. 
\item For some $\phi(x)$ cohomologous to
$\Phi(x,y)$  holds $P(\phi(x)) < \infty$. 
\item Norm $\norm{L_{\phi}(1)}_{\infty}$ is finite.
\ee
then \ref{maindecay} holds.
\end{Prop}
Proposition \ref{sufficient} gives sufficient conditions
for exponential decay of correlations for H\"older (in particular smooth)
functions restricted to the attractor. 
\ee
\section{ Proof of the exponential decay of correlations}
We check properties $(a)$ - $(d)$.
\be
\item
Recall that in our model we consider the   partition of the square into full height rectangles $E_i$.\\
Our shift is Bernoulli, all rows (and columns)  are the same row of $1$-s, so 
it is topologically mixing and property
$(a)$ is satisfied.
\item Property $(b)$ follows from \ref{var1}. 
\item Next we prove property $(c)$. As in the case of attractors for Axiom A  systems
we prove 
\begin{Prop} \label {Pequals0} 
For $\phi(x) = \phi^u(x)$ topological pressure  $P(\phi^u(x)) $ equals zero.
\end{Prop}
{\bf Proof}. \\
We fix some symbol $a$,  respective rectangle $E_a$, and $W^u_{0a}=W^u_0 \cap E_a$. 
When evaluating  $Z_n(\phi,a) $  in \ref{statsum} we consider respective sum
over all periodic orbits of period $n$ starting in $E_a$. \\
Each cylinder set ${\Ean}$ contains one periodic orbit of period $n$. 
When evaluating  $\phi^u(x)$ we use formula \ref{cohomformula}. We can evaluate that expression 
at a point $z$ 
of intersection between the stable manifold of a periodic point in ${\Ean}$ and $W^u_{0a}$.
Then each term in \ref{statsum}   is a product
of two expressions. \\ 
The first expression equals $\frac{1}{D^uF^n(z)}$, which  
coincides up to a uniformly bounded factor
with the length of $W^u(z,{\Ean})$. \\
The second expression $e^{ u(z) - u(F(z)}$ is uniformly bounded away from zero and infinity. \\
Therefore up to  a uniformly bounded factor the sum \ref{statsum} equals to the length of  $W^u_{0a}$.
That implies $ P(\phi)=0$. \ \ Q.E.D. \\
\item
Property $(d)$  is an easier version of the above proposition.
For any $x$  the expression $L_{\phi}(1(x))$ equals up to a uniform constant
\[\sum_i \norm{W^u_0(x_i,E_i)} \]
where $x_i \in E_i$ is the point of intersection between the stable manifold of the respective 
preimage of $x$ and $W^u_0$. Therefore $\norm{L_{\phi}(1)}_{\infty}< \infty$
 is bounded by $c\norm{W^u_0}$.
\ee
So all properties of Proposition \ref{sufficient} are satisfied, and we get exponential decay of
correlations for one-sided shift. As in \cite{Bowen} it implies
exponential decay of
correlations for two-sided shift and therefore for H\"older functions on $Q$. That
 proves Theorem \ref{Main}. \\
\begin{Rem} \ref{maindecay}  also implies
the central limit theorem  for  for H{\"o}lder functions on $Q$ .
\end{Rem}

 Let us denote $\mu_S$ the invariant measure on $\Lambda$ constructed in \cite{Jak-New-1}, \cite{Jak-New-2}
following Sinai method, and let $\mu_{RB}$ be the invariant measure on $\Lambda$
constructed above  following Ruelle-Bowen method, see \cite{Ruelle-1}, \cite{Bowen}.
Let $\mu_1$ be the projection of $\mu_S$ onto one-sided sequences, and let $\mu$
be the measure on one-sided sequences constructed above by Ruelle-Bowen method.
In both  constructions measures of cylinder sets $[i_0i_1 \ldots i_{n-1}]$ of any rank
equal up to a uniform constant to the length of the crossections\
of ${\En}$ by  $W^u_0$.
So $\mu_1$ and $\mu$ are equivalent and therefore they coincide. That is a particular case
of the characterization of Gibbs measures proved in \cite{Sarig-1}. \\
As in the classical case that implies
\begin{Cor}
Measures $\mu_S$ and $\mu_{RB}$ coincide. 
\end{Cor}

{\bf Acknowledgements. } I want to thank Sheldon Newhouse, David Ruelle and 
Omri Sarig for  useful discussions during the preparation of this paper.

%%%%%%%%%%%%%%%%%%%%%%%%%%%%%%%%
%%%%%%%%%%%%%%%%%%%%%%%%%%%%%%%%%%

\end{document}